# A Sparse Bayesian Estimation Framework for Conditioning Prior Geologic Models to Nonlinear Flow Measurements

Lianlin Li and Behnam Jafarpour

*Department of Petroleum, Texas A&M University, 77845, USA*

## ABSTRACT

We present a Bayesian framework for reconstruction of subsurface hydraulic properties from nonlinear dynamic flow data by imposing sparsity on the distribution of the solution coefficients in a compression transform domain. Sparse representation of the subsurface flow properties in a compression transform basis lends itself to a natural regularization approach, i.e. sparsity regularization, which has recently been exploited in solving ill-posed nonlinear inverse problems that frequently encountered in subsurface flow and transport modeling. The Bayesian estimation approach allows for a probabilistic treatment of the sparse reconstruction problem by enforcing sparsity through Laplace priors on the distribution of the solution in the sparsifying transform basis. The methodology has its roots in machine learning and recently introduced relevance vector machine algorithm for *linear* inverse problems. We extend the application of this approach to *nonlinear* subsurface inverse problems where solution sparsity in a discrete cosine transform is assumed. The probabilistic fulfillment of solution sparsity, as opposed to deterministic regularization, avoids the nuisance of specifying *a priori* regularization parameter and allows for quantification of the estimation and prediction uncertainty. Several numerical experiments from subsurface multiphase flow and transport application are conducted to illustrate the performance of proposed method and compare it with the regular Bayesian estimation approaches that do not impose solution sparsity. While the examples are derived from a geophysical application, the proposed framework can be applied to *nonlinear* inverse problems in other fields such as medical imaging and electromagnetic inverse problem.

# I. INTRODUCTION

The fluid flow and transport equations in porous media that are derived from Darcy's law and mass conservation principle are widely used to quantify and predict fluid displacement behavior in the subsurface environment [1-4, 23, 36, 37]. A critical component of the flow prediction is model calibration, which refers to identification of key model input parameters from observed flow quantities to improve the prediction accuracy. Hydraulic rock properties such as permeability and porosity play an essential role in controlling the fluid displacement patterns and are therefore the main targets of the model calibration study. The general framework for performing model calibration involves minimization of a prescribed cost function that penalizes the misfit between measured and predicted data and deviations from a best-known prior model by adapting the set of sensitive parameters (e.g. permeability). In general, the estimation of heterogeneous hydraulic rock properties such as intrinsic permeability from dynamic flow measurements leads to an underdetermined nonlinear inverse problem [10, 11, 13, 14, 18-21, 26-28, 36, 37]. That is, the number of unknown parameters is significantly larger than available data. Consequently, several non-unique solutions can exist that explain the available data equally well but fail to predict the future flow behavior correctly [10, 11, 13, 14, 18-21, 26-28, 36, 37].

To improve the non-uniqueness of the inverse modeling solution, additional sources of information, mainly in the form of structural prior information, are incorporated into the solution framework to remove irrelevant solutions.. While spatial regularization techniques [5, 10, 11, 18-21, 26, 31-33] are commonly used to improve the solution of such ill-posed inverse problems, in some cases model representation in an appropriate transformed domain such as Fourier can lead to more suitable and natural regularization formulations without imposing simplified prior structural assumptions (i.e. smoothness) on the solution [5, 10, 11, 18-21, 26, 31-33]. Image compression transform bases such as Fourier [6] and Wavelet [22] have been widely applied to parsimoniously represent spatially correlated images/volumes with a sparse set of coefficients. Because of the continuity (i.e. spatial correlation) in geologic facies, subsurface features are amenable to sparse approximations in a compression basis. The sparsifying nature of these transforms provides a general regularization framework in the transform domain, i.e. sparsity regularization, which has recently been used to regularize ill-posed inverse problems [7-9, 12, 16, 20, 21, 25, 29, 32, 35]. Owing to the sparse approximation of continuous geologic facies in an appropriate transform basis, a quest for sparse solutions that match available inversion data is a natural consequence without placing too restrictive and simplifying structural assumptions on the spatial representation of the solution. The expected outcome of this solution framework is to identify, as warranted by the available data, the main *features* in the solution without imposing strict *a priori* structural assumptions. As discussed in the next section, the search for a sparse solution is facilitated by recent developments in sparse signal reconstruction known as the *compressed sensing* paradigm [7-9, 12, 16, 29].

An alternative approach to deterministic regularization methods is the Bayesian estimation framework that allows for a probabilistic treatment of the prior knowledge and observations to characterize the inversion solution and its statistical distribution. The Bayesian estimation has recently been used to solve linear inverse problem with sparsity promoting constraints in machine learning, leading to the sparse Bayesian learning algorithm [29, 34] We extend the application of the sparse Bayesian estimation approach *nonlinear dynamic* inverse problem where reconstruction of permeability fields from nonlinear dynamic data is considered. We achieve this through iterative solution of the linearized form of the nonlinear likelihood function. For the inverse problem studied, the measurements include pressures and saturations at well (sink/source) locations. These measurements have different scales and statistical properties, which we model using the Gaussian mixture prior, where each measurement has its own noise variance.

The reminder of this paper is arranged as follows. In section II, the inverse problem and sparsity-promoting regularization of it are briefly described. A detailed presentation of the sparse Bayesian framework is presented in section III followed by the nonlinear application derived in this paper. Section IV is dedicated to numerical experiments to illustrate the nonlinear inversion performance of developed approach before the final remarks and conclusions are presented in Section V.

## II. PROBLEM STATEMENT

Identification of flow-related heterogeneous subsurface properties such as the spatial distribution intrinsic rock permeability from dynamic flow data naturally leads to an underdetermined inverse problem [10, 11, 13, 14, 18-21, 26-28, 36, 37] because the number of unknown variables significantly exceeds the number of available measurements. Therefore, one is always confronted by the possibility of many non-unique solutions that adequately describe the existing measurements. To alleviate solution non-uniqueness, it is common to reduce the number of parameters (i.e. parameterization) or to introduce additional assumptions about the solution to regularize the problem (i.e. regularization) or both [10, 11, 13, 14, 18-21, 26-28, 36, 37]. In either case, to increase the likelihood of finding only relevant solutions, one has to incorporate, often based on the physics of the problem, additional prior assumptions such as smoothness or flatness [5, 31, 33, 36]. Several authors [31] have extensively discussed approaches to deal with underdetermined inverse problems, especially within geophysical inverse modeling and subsurface characterization context that readers can refer to for details. Here, we begin by introducing a typical regularization formulation of the nonlinear inverse problem in subsurface flow modeling before introducing the Sparse Bayesian Learning (SBL) [29, 34, 35] framework.

In a typical nonlinear inverse problem that arise in multiphase flow in porous media application the solution $\mathbf{m} \in \mathbb{R}^N$ of an unknown rock property (e.g. permeability) is obtained by minimizing a cost function $J(\mathbf{m})$ that consists of a data misfit penalty term and usually one or more regularization term(s) that penalize deviation from a specified solution structure [10, 11, 13, 14, 18-21, 26-28, 36, 37]. Without the regularization term, the cost-functional to be minimized over $\mathbf{m}$ is of the form $J(\mathbf{m}) = \|\mathbf{y} - \mathbf{g}(\mathbf{m})\|_2^2$, where $\mathbf{y} = \begin{bmatrix} \mathbf{y}_S \\ \mathbf{y}_P \end{bmatrix} \in \mathbb{R}^M$ is the vector of measured data with $\mathbf{y}_P$ and $\mathbf{y}_S$ representing point measurements of pressure and saturation and $\mathbf{y}^{sim} = \mathbf{g}(\mathbf{m}) = \begin{bmatrix} \mathbf{g}_S(\mathbf{m}) \\ \mathbf{g}_P(\mathbf{m}) \end{bmatrix} \in \mathbb{R}^M$ is the vector of corresponding predicted (simulated) measurements.

In the examples considered in this paper, the simulated data are nonlinearly related to model parameters $\mathbf{m}$, hence the resulting objective function $J$ becomes a nonlinear function of unknown parameters. A simple regularization approach is to augment the above data misfit cost function with term that encourages simple solutions with minimal energy, i.e. $\|\mathbf{m}\|_2^2$. The resulting cost-function is

$$J(\mathbf{m}) = \|\mathbf{y} - \mathbf{g}(\mathbf{m})\|_2^2 + \delta \|\mathbf{\Phi}\mathbf{m}\|_2^2 \qquad (1)$$

where $\delta$ is a regularization parameter that controls the relative weight given to the regularization term. Other commonly used spatial regularization techniques that promote smoothness or flatness of the solution (known as Tikhonov regularization [33]) have a similar form to (1), except for an operator $\mathbf{L}$ (such as first or second order spatial derivative) that is applied to $\mathbf{m}$ to define the desired structure (i.e. the norm of $\mathbf{Lm}$ is used as the regularization term).

An important issue in practical implementation of the regularized inverse problem above is the lack of knowledge about the regularization parameter. In general, a reasonable choice for the regularization parameter can only be obtained either through a priori knowledge about the solution or through guidelines specified in methods such as L-Curve or the generalized cross validation (GCV) technique [31, 33] that require extensive numerical experimentation. Aside from this nuisance, regularization provides a simple approach to find *a local* deterministic solution to the ill-posed inverse problem. The above spatial regularization approach can also be applied in a transform domain if the solution is expected to have quantifiable transform-domain structural attributes. Jafarpour et al. [20, 21] considered *sparsity* of the rock permeability solution in the discrete cosine transform (DCT), a Fourier-related transform, domain as a general structural attribute of the solution that they used to formulate an effective transform-domain regularization technique. The formulation is inspired from recent advances in sparse reconstruction literature [9, 12, 25, 32] and the spatial continuity of rock properties that translates into spatial correlation, which lend themselves to parsimonious/compressed transform-domain representations. Examples of decorralting or compressive transforms are the preconstructed Fourier or Wavelet such as the DCT or DWT (discrete wavelet transform) bases that are used in image and video compression [6, 22].

To illustrate this point, Fig. 1a and 1b show a typical permeability distribution with its corresponding DWT and DCT coefficients, respectively. An approximate reconstruction of the permeability field in each figure with 1%, 2% and 5% largest coefficients are also shown in the second through fourth column of each figure. This simple example suggests that a small fraction of the transformed coefficients can capture the most salient features in the original image. The DCT and DWT bases are only two of many other transform basis with excellent compression property that can be used to sparsify correlated spatial images such as the permeability example shown. In this paper, we assume that a suitable sparsifying basis $\boldsymbol{\Phi}$ exists to compress the spatial representation of the inverse problem solution and formulate a Sparse Bayesian Learning (SBL) [34] approach for finding the sparse representation from *nonlinear* dynamic data. We begin our formulation by briefly presenting the notations.

We begin the SBL formulation by reviewing the sparsity-promoting regularization problem in a transform domain, whose original formulation is expressed as

$$\min_{\mathbf{m}} \|\boldsymbol{\Phi}\mathbf{m}\|_0 \tag{P0}$$

$$s.t. \quad \mathbf{y} = \mathbf{g}(\mathbf{m}) \quad \text{or} \quad \|\mathbf{y} - \mathbf{g}(\mathbf{m})\|_2^2 \leq \sigma$$

where $\sigma$ is the noisy energy and minimization of the $l_0$-quasinorm, which counts the number of nonzero elements of its argument, encourages solutions with minimum support. The exact solution to (P0) requires a combinatorial search over across all possible sparse sets of the basis, which is known to be NP-hard complex [7, 8, 16]. A vast and growing literature exists on approximate solutions to (P0). The practical solution techniques can be classified as methods that apply a *convex relaxation* [7-9, 16] to (P0) and solve the corresponding convex optimization problem and *greedy algorithms* known as matching pursuit [16]. Donoho et al. [1-9, 16] showed that when $g(m)$ is linear solution guarantees exist through $l_1$-norm convexification of (P0) for sufficiently sparse solutions and with adequate measurements. The resulting formulation can be expressed as

$$\min_{\mathbf{m}} \|\mathbf{\Phi m}\|_1 \quad (P1)$$

$$s.t. \quad \mathbf{y} = \mathbf{g}(\mathbf{m}) \quad or \quad \|\mathbf{y} - \mathbf{g}(\mathbf{m})\|_2^2 \leq \sigma$$

The regularized form of (P1) is

$$\min_{\mathbf{m}} L = \|\mathbf{y} - \mathbf{g}(\mathbf{m})\|_2^2 + \delta \|\mathbf{\Phi m}\|_1 \quad (P2)$$

where $\|\mathbf{x}\|_1 = \sum_{i=1}^{N} |\mathbf{x}_i|$ denotes the $l_1$-norm of vector $\mathbf{x}$. The resulting transform-domain regularized inverse problem can be solved to find a sparse approximation to the original parameters in space. Application of this formulation to reconstruct permeability fields from dynamic flow measurement has been discussed in [20, 21]. In the next section, we present a more general Bayesian formulation of the sparse reconstruction problem known as Sparse Bayesian Learning.

### III. SPARSE BAYESIAN MODELLING FOR HISTORY MATCHING

Problems (P1) and (P2) can be solved within the deterministic regularization approach described above. However, the regularized formulation of the problem suffers from a number of issues including specification of regularization parameter and the deterministic nature of the formulation. The lack of sufficient data to uniquely constrain a prior model calls for a solution approach that can quantify the expected error in the estimation or provides multiple likely solutions. Furthermore, the uncertainty in the prior model input parameters and the noise in the measurements are better dealt with in a probabilistic data integration framework. The Bayesian [17, 18, 29, 34, 35] framework provides an elegant formulation for combining uncertain sources of information through their probabilistic representations. The Bayesian formulation provides distinct advantages over deterministic framework by allowing for probabilistic predictions, automatic incorporation and estimation of uncertain model parameters, and estimation of the reconstruction uncertainty.

One approach to formulate a sparsity-promoting Bayesian framework is to specify independent Laplace priors on the unknown parameters in the transform domain. We note that the DCT coefficients of natural images are known to follow a Laplace distribution. Furthermore, because of the decorrelating property of the DCT basis, the transformed coefficients are approximately uncorrelated, hence individual Laplace distributions are applied to each transformed coefficients. Denoting the full transform domain representation of parameter (permeability) $\mathbf{m}$ by the vector $\boldsymbol{\alpha}$, we have $\boldsymbol{\alpha} = \boldsymbol{\Phi}\mathbf{m}$, where we take the matrix $\boldsymbol{\Phi}$ to have as its columns the full DCT basis vectors. However, if a Laplace prior is assumed for each individual DCT coefficients, the resulting problem is not easily tractable because Laplace priors are not conjugate functions to Gaussian likelihoods (commonly assumed distribution for likelihood functions), i.e. under Gaussian likelihood, the resulting posterior distribution will not be a Laplace distribution [34]. To circumvent this issue, we impose the Laplace prior on the transformed coefficients using a hierarchical approach or Gaussian mixtures as proposed in [29, 34]; furthermore, we assume independent Gaussian distributions for each measurement $y_i$. The details of the SBL formulation are presented next.

**III.1 Models of permeability and observations**

In the "Bayesian" land, all unknown/uncertain variables, including observation noise, the parameters ( permeability in this paper) and other derived hyperparameters (will be defined soon), are treated as stochastic quantities with assigned probability distributions. We present each of the involved variables and their specified distribution in this section.

Firstly, we model each observation noise $\mathbf{n} = \mathbf{y} - \mathbf{g}(\mathbf{m})$ as an independent Gaussian process with zero mean and unknown variance, that is,

$$\Pr(\mathbf{y} \mid \mathbf{m}, \mathbf{B}) = \mathcal{N}\left(\mathbf{y} \mid \mathbf{g}(\mathbf{m}), \mathbf{B}^{-1}\right) \tag{2}$$

where $\mathbf{B} = diag\{\boldsymbol{\beta}_i, i = 1, 2, ..., M\}$. We treat the observation variances $\mathbf{B}$ as unknown hyperparameters in this formulation. Since the Gamma distribution is the conjugate prior for inverse variance of *Gaussian* distributions, the analysis is greatly simplifed by placing the Gamma prior on $\boldsymbol{\beta}_i$, i.e.

$$\Pr(\boldsymbol{\beta}_i \mid a_\beta, b_\beta) = Gamma(\boldsymbol{\beta}_i \mid a_\beta, b_\beta) \quad, \quad i = 1, 2, ..., M \tag{3}$$

with the *Gamma* distribution is defined as

$$Gamma(x\mid a,b) := \frac{b^a}{\Gamma(a)} x^{a-1} \exp(-bx), \quad x \geq 0 \qquad (4)$$

To obtain large $\boldsymbol{\beta}_i$ values, we may choose $a_\beta > 1$ and $b_\beta$ to be very small. In this paper, $a_\beta = 1.0 + \varepsilon$ and $b_\beta = \varepsilon$ are specified, where $\varepsilon > 0$ is a very small value, and denoted as $a_\beta := 1_+$ and $b_\beta := 0_+$.

Moving next to the prior permeability model, to realize $l_1$ regularization formulation in (P1) within the Bayesian estimation framework, the Laplace prior should be imposed on each coefficient $\boldsymbol{\alpha}_i$, that is,

$$\Pr(\boldsymbol{\alpha}_i \mid \boldsymbol{\lambda}_i) = \frac{\lambda_i}{2} \exp\left(-\frac{\lambda_i}{2}|\boldsymbol{\alpha}_i|\right), \quad i = 1, 2, \ldots, N \qquad (5)$$

However, since this distribution is not conjugate to the conditional distribution in (2), it does not allow for a tractable Bayesian analysis. Hierarchical priors or Gaussian mixtures are employed be to describe $\boldsymbol{\alpha}_i$ [34]. As the first stage of a hierarchical model, independent Gaussian priors are specified for each coefficient $\boldsymbol{\alpha}_i$,

$$\Pr(\boldsymbol{\alpha}_i \mid \boldsymbol{\gamma}_i) = \mathcal{N}(\boldsymbol{\alpha}_i \mid 0, \boldsymbol{\gamma}_i) \qquad (6)$$

To make the priors on the hyperparameters non-informative (i.e. flat), we may choose $\boldsymbol{\gamma}_i$ to be very small. Using small values for the hyperparameter allows the posterior probability of $\boldsymbol{\alpha}_i$s to concentrate at very large values, resulting in many of the DCT basis to have zero coefficients and effectively leading to sparsity by removing irrelevant basis elements. In the second stage of the hierarchy, the Gamma distribution with $\left(1, \frac{\lambda_i}{2}\right)$ are assigned to $\boldsymbol{\gamma}_i$, in particular,

$$\Pr(\gamma_i \mid \lambda_i) = Gamma\left(\gamma_i \mid 1, \frac{\lambda_i}{2}\right) \tag{7}$$

Combing (7) and (6) results in (5) due to $\Pr(\alpha_i \mid \lambda_i) = \int \Pr(\alpha_i \mid \gamma_i) \Pr(\gamma_i \mid \lambda_i) d\gamma_i$, which imposes the sparsity-promoting constraint on $\boldsymbol{\alpha}$ via $l_1$-norm regularization.

**III.2 Bayesian inference**

Having defined the prior distribution for the hyperparameters, Bayesian inference can now be carried out to compute the posterior distribution over all unknowns, including $\{\boldsymbol{\alpha}_i, i = 1, 2, ..., N\}$, $\{\boldsymbol{\beta}_i, i = 1, 2, ..., M\}$ and the derived hyperparameters given the data. The $\boldsymbol{\alpha}$ is estimated from the posterior of $\boldsymbol{\alpha}$ over observation $y$ and hyperparameters $\Lambda := diag\{\gamma_i^{-1}, i = 1, 2, ..., N\}$ and $\mathbf{B}$, in particular,

$$\boldsymbol{\alpha} = \arg\max_x \Pr(\mathbf{x} \mid \mathbf{y}, \mathbf{B}, \Lambda) \tag{8}$$

where

$$\begin{aligned}\Pr(\mathbf{x} \mid \mathbf{y}, \mathbf{B}, \Lambda) &\propto \Pr(\mathbf{y} \mid \mathbf{x}, \mathbf{B}) \Pr(\mathbf{x} \mid \Lambda) \\ &\propto \exp\left(-\frac{1}{2}\|\mathbf{y} - \mathbf{g}(\mathbf{m})\|_\mathbf{B}^2 - \frac{1}{2}\|\Phi \mathbf{m}\|_\Lambda^2\right)\end{aligned} \tag{9}$$

Iterative approximation is introduced to deal with the nonlinearity of $\mathbf{g}(\mathbf{m})$. From the first-order Taylor approximation of $\mathbf{g}(\mathbf{m}^{(n+1)})$ at iteration $n+1$, we have

$$\|\mathbf{y} - \mathbf{g}(\mathbf{m}^{(n+1)})\|_B^2 \approx \|\mathbf{y}^{(n)} - \mathbf{G}^{(n)} \mathbf{m}^{(n+1)}\|_B^2 \tag{10}$$

with $\mathbf{y}^{(n)} = \mathbf{y} - \mathbf{g}(\mathbf{m}^{(n)}) + \mathbf{G}^{(n)} \mathbf{m}^{(n)}$, where $\mathbf{G}^{(n)} = \begin{bmatrix} \mathbf{G}_S^{(n)} \\ \mathbf{G}_P^{(n)} \end{bmatrix}$ is the sensitivity matrix that can be efficiently obtained from the solution of the forward flow simulation. Correspondingly, (8) is expressed as

$$\boldsymbol{\alpha}^{(n+1)} = \arg\max_{\alpha} \Pr\left(\boldsymbol{\alpha} \mid \mathbf{y}^{(n)}, \mathbf{B}^{(n+1)}, \boldsymbol{\Lambda}^{(n+1)}, \mathbf{m}^{(n)}\right) \tag{11}$$

where $\mathbf{B}^{(n+1)}$ and $\boldsymbol{\Lambda}^{(n+1)}$ are the noise inverse covariance and signal covariance, respectively, at the ($n+1$)-th iteration. After simple algebraic manipulations, we obtain

$$\begin{aligned}\Pr\left(\boldsymbol{\alpha} \mid \mathbf{y}, \mathbf{B}^{(n+1)}, \boldsymbol{\Lambda}^{(n+1)}, \mathbf{m}^{(n)}\right) &\propto \exp\left(-\frac{1}{2}\left\|\mathbf{y}^{(n)} - \mathbf{G}^{(n)}\mathbf{m}\right\|^2_{\mathbf{B}^{(n+1)}} - \frac{1}{2}\left\|\boldsymbol{\Phi}\mathbf{m}\right\|^2_{\boldsymbol{\Lambda}^{(n+1)}}\right) \\ &\propto \left(\frac{1}{2\pi}\right)^{N/2} |\boldsymbol{\Sigma}|^{-0.5} \exp\left(-\frac{1}{2}\|\boldsymbol{\alpha}-\boldsymbol{\mu}\|^2_{\boldsymbol{\Sigma}^{-1}}\right) \times \left(\frac{1}{2\pi}\right)^{M/2} |\mathbf{C}|^{-0.5} \exp\left(-\frac{1}{2}\left\|\mathbf{y}^{(n)}\right\|^2_{\mathbf{C}^{-1}}\right)\end{aligned} \tag{12}$$

where

$$\boldsymbol{\mu} = \boldsymbol{\Sigma}\tilde{\mathbf{G}}^T \mathbf{B}^{(n+1)}\mathbf{y}, \tag{13}$$

$$\boldsymbol{\Sigma} = \left(\boldsymbol{\Lambda}^{(n+1)} + \tilde{\mathbf{G}}^T \mathbf{B}^{(n+1)}\tilde{\mathbf{G}}\right)^{-1}, \tag{14}$$

$$\mathbf{C} = \left[\mathbf{B}^{(n+1)}\right]^{-1} + \tilde{\mathbf{G}}\left[\boldsymbol{\Lambda}^{(n+1)}\right]^{-1}\tilde{\mathbf{G}}^T = \left(\mathbf{B}^{(n+1)} - \mathbf{B}^{(n+1)}\tilde{\mathbf{G}}\boldsymbol{\Sigma}\tilde{\mathbf{G}}^T\mathbf{B}^{(n+1)}\right)^{-1} \tag{15}$$

$$\tilde{\mathbf{G}} = \mathbf{G}^{(n)}\boldsymbol{\Phi}^T \tag{16}$$

In the above equations, we have assumed, without loss of generality, that the basis $\boldsymbol{\Phi}$ is orthogonal. The posterior covariance and mean of $\alpha$ at iteration $n+1$ are expressed, respectively, as

$$\boldsymbol{\Sigma} = \left(\boldsymbol{\Lambda}^{(n+1)} + \tilde{\mathbf{G}}^T \mathbf{B}^{(n+1)}\tilde{\mathbf{G}}\right)^{-1} \tag{17}$$

and

$$\boldsymbol{\alpha}^{n+1} = \boldsymbol{\Sigma}\tilde{\mathbf{G}}^T \mathbf{B}^{(n+1)}\mathbf{y}^{(n)} \tag{18}$$

The Bayesian estimation of hyperparaperters via the so-called *evidence procedure* (type-II maximum likelihood approach), is considered next. In particular, Bayesian inference is performed based on the following decomposition at ($n+1$)th iteration,

$$\left(\hat{\mathbf{\Lambda}}^{(n+1)}, \hat{\mathbf{B}}^{(n+1)}, \hat{\boldsymbol{\lambda}}^{(n+1)}\right) = \arg\max_{\left\{\mathbf{\Lambda}^{(n+1)}, \mathbf{B}^{(n+1)}, \boldsymbol{\lambda}^{(n+1)}\right\}} \Pr\left(\mathbf{\Lambda}^{(n+1)}, \mathbf{B}^{(n+1)}, \boldsymbol{\lambda}^{(n+1)} \mid \mathbf{y}^{(n)}\right) \quad (19)$$

where

$$\begin{aligned}
\Pr&\left(\mathbf{\Lambda}^{(n+1)}, \mathbf{B}^{(n+1)}, \boldsymbol{\lambda}^{(n+1)} \mid \mathbf{y}^{(n)}\right) \\
&\propto \Pr\left(\mathbf{y}^{(n)} \mid \mathbf{\Lambda}^{(n+1)}, \mathbf{B}^{(n+1)}\right) \Pr\left(\mathbf{\Lambda}^{(n+1)}, \mathbf{B}^{(n+1)}, \boldsymbol{\lambda}^{(n+1)}\right) \\
&= \Pr\left(\mathbf{y}^{(n)} \mid \mathbf{\Lambda}^{(n+1)}, \mathbf{B}^{(n+1)}\right) \Pr\left(\mathbf{\Lambda}^{(n+1)} \mid \boldsymbol{\lambda}^{(n+1)}\right) \Pr\left(\mathbf{B}^{(n+1)}\right) \\
&= \left(\frac{1}{2\pi}\right)^{M/2} |\mathbf{C}|^{-0.5} \exp\left(-\frac{1}{2}\|\mathbf{y}^{(n)}\|^2_{\mathbf{C}^{-1}}\right) \times \left[\prod_{i=1}^{N} Gamma\left(\gamma_i \mid 1, \frac{\lambda_i}{2}\right)\right] \\
&\quad \times \left[\prod_{i=1}^{M} Gamma\left(\beta_i \mid 1_+, 0_+\right)\right]
\end{aligned} \quad (20)$$

Maximization of the logarithm of (20) leads to maximizing

$$\left(\hat{\mathbf{\Lambda}}^{(n+1)}, \hat{\mathbf{B}}^{(n+1)}, \hat{\boldsymbol{\lambda}}^{(n+1)}\right) = \arg\max_{\left\{\mathbf{\Lambda}^{(n+1)}, \mathbf{B}^{(n+1)}, \boldsymbol{\lambda}^{(n+1)}\right\}} \log\left[\Pr\left(\mathbf{\Lambda}^{(n+1)}, \mathbf{B}^{(n+1)}, \boldsymbol{\lambda}^{(n+1)} \mid \mathbf{y}^{(n)}\right)\right] \quad (21)$$

where

$$\begin{aligned}
\mathcal{L} &:= \log \Pr\left(\mathbf{\Lambda}^{(n+1)}, \mathbf{B}^{(n+1)}, \boldsymbol{\lambda}^{(n+1)} \mid \mathbf{y}^{(n)}\right) \\
&\propto -\frac{1}{2}\log|\mathbf{C}| - \frac{1}{2}\|\mathbf{y}^{(n)}\|^2_{\mathbf{C}^{-1}} + \sum_{i=1}^{N} \log Gamma\left(\gamma_i \mid 1, \frac{\lambda_i}{2}\right) + \sum_{i=1}^{M} \log Gamma\left(\beta_i \mid 1_+, 0_+\right)
\end{aligned} \quad (22)$$

Updating other hyperparameters is carried out by taking the derivative of (22) with respect to each hyperparameter and setting it equal to zero. To do this, we should find the derivative of $-\frac{1}{2}\log|\mathbf{C}| - \frac{1}{2}\|\mathbf{y}^{(n)}\|^2_{\mathbf{C}^{-1}}$ with respect to $\beta_i$ and $\gamma_i$. Upon using the Woodbury identity [34], we get

$$\begin{aligned}
\|\mathbf{y}^{(n)}\|^2_{\mathbf{C}^{-1}} &= \left[\mathbf{y}^{(n)}\right]^T \left[\mathbf{B}^{(n+1)} - \mathbf{B}^{(n+1)}\tilde{\mathbf{G}}\boldsymbol{\Sigma}\tilde{\mathbf{G}}^T\mathbf{B}^{(n+1)}\right] \mathbf{y}^{(n)} \\
&= \left[\mathbf{y}^{(n)}\right]^T \mathbf{B}^{(n+1)} \left(\mathbf{y}^{(n)} - \tilde{\mathbf{G}}\boldsymbol{\mu}\right) \\
&= \left\|\mathbf{y}^{(n)} - \tilde{\mathbf{G}}\boldsymbol{\mu}\right\|^2_{\mathbf{B}^{(n+1)}} + \boldsymbol{\mu}^T\tilde{\mathbf{G}}^T\mathbf{B}^{(n+1)}\left(\mathbf{y}^{(n)} - \tilde{\mathbf{G}}\boldsymbol{\mu}\right) \\
&= \left\|\mathbf{y}^{(n)} - \tilde{\mathbf{G}}\boldsymbol{\mu}\right\|^2_{\mathbf{B}^{(n+1)}} + \boldsymbol{\mu}^T\mathbf{\Lambda}^{(n+1)}\boldsymbol{\mu}
\end{aligned} \quad (23)$$

leading to

$$\frac{\partial \|\mathbf{y}^{(n)}\|^2_{\mathbf{C}^{-1}}}{\partial \boldsymbol{\beta}_i} = \left(\mathbf{y}^{(n)} - \tilde{\mathbf{G}}\boldsymbol{\mu}\right)^2_i \tag{24}$$

and

$$\frac{\partial \|\mathbf{y}^{(n)}\|^2_{\mathbf{C}^{-1}}}{\partial \boldsymbol{\gamma}_i} = \boldsymbol{\mu}_i^2 \tag{25}$$

On the other hand, the gradient of $\log|\mathbf{C}|$ with respect to $\boldsymbol{\beta}_i$ is

$$\frac{\partial \log|\mathbf{C}|}{\partial \boldsymbol{\beta}_i} = \frac{1}{|\mathbf{C}|}\frac{\partial |\mathbf{C}|}{\partial \boldsymbol{\beta}_i} = \frac{1}{|\mathbf{C}|}\frac{\partial \left|\left[\mathbf{B}^{(n+1)}\right]^{-1} + \tilde{\mathbf{G}}\left[\mathbf{\Lambda}^{(n+1)}\right]^{-1}\tilde{\mathbf{G}}^T\right|}{\partial \boldsymbol{\beta}_i} \tag{26}$$
$$= -\boldsymbol{\beta}_i^{-2}\left[\mathbf{C}^{-1}\right]_{ii}$$

where $\left[\mathbf{C}^{-1}\right]_{ii}$ is the $(i,i)$-entry of $\mathbf{C}^{-1}$, which can be simplified into

$$\left[\mathbf{C}^{-1}\right]_{ii} = \left(\mathbf{B}^{(n+1)} - \mathbf{B}^{(n+1)}\tilde{\mathbf{G}}\boldsymbol{\Sigma}\tilde{\mathbf{G}}^T\mathbf{B}^{(n+1)}\right)_{ii} = \boldsymbol{\beta}_i - \boldsymbol{\beta}_i^2\left(\tilde{\mathbf{G}}\boldsymbol{\Sigma}\tilde{\mathbf{G}}^T\right)_{ii} \tag{27}$$

Combing (22), (24), (26) and (27) yields the estimate for

$$\mathbf{B}^{(n+1)} = diag\{\boldsymbol{\beta}_i, i=1,2,...,M\} \text{ with } \boldsymbol{\beta}_i = \frac{1}{\left(\tilde{\mathbf{G}}\boldsymbol{\Sigma}\tilde{\mathbf{G}}^T\right)_{ii} + \left(\mathbf{y}^{(n)} - \tilde{\mathbf{G}}\boldsymbol{\mu}\right)^2_i} \tag{28}$$

Several remarks regarding (28) are in order

(i) To avoid possible singularity due to small values of the denominator in (28), we introduce a small value $\varepsilon > 0$, to obtain the modified formulation for $\boldsymbol{\beta}_i$

$$\boldsymbol{\beta}_i = \frac{1}{\left(\tilde{\mathbf{G}}\boldsymbol{\Sigma}\tilde{\mathbf{G}}^T\right)_{ii} + \left(\mathbf{y}^{(n)} - \tilde{\mathbf{G}}\boldsymbol{\mu}\right)^2_i + \varepsilon} \tag{29}$$

From our numerical experiences, the final reconstruction is sensitivity to this choice. In this paper, $\varepsilon = 10^{-6}$ is specified.

(ii) It is easy to check that $\boldsymbol{\beta}_i$ depends on the misfit $\left(\mathbf{y}^{(n)} - \tilde{\mathbf{G}}\boldsymbol{\mu}\right)^2_i$, which is used to balance the importance of measurements. Inspired by this relation, we propose another version of the above

Bayesian inference formulation; by introducing $\eta_S = \dfrac{1}{\left\|\mathbf{y}_S - \mathbf{g}_S\left(\mathbf{m}^{(n)}\right)\right\|_2}$ and $\eta_P = \dfrac{1}{\left\|\mathbf{y}_P - \mathbf{g}_P\left(\mathbf{m}^{(n)}\right)\right\|_2}$, the misfit and the sensitivity matrix at (n+1)th iteration are modified as one in the weighting form, in particular,

$$\mathbf{G}^{(n)} \leftarrow \mathbf{G}^{(n)} = \begin{bmatrix} \eta_S \mathbf{G}_S^{(n)} \\ \eta_P \mathbf{G}_P^{(n)} \end{bmatrix}, \tag{30}$$

$$\mathbf{y} = \begin{bmatrix} \mathbf{y}_S \\ \mathbf{y}_P \end{bmatrix} \leftarrow y = \begin{bmatrix} \eta_S \mathbf{y}_S \\ \eta_P \mathbf{y}_P \end{bmatrix}, \tag{31}$$

and

$$\mathbf{g}(\mathbf{m}) = \begin{bmatrix} \mathbf{g}_S(\mathbf{m}) \\ \mathbf{g}_P(\mathbf{m}) \end{bmatrix} \leftarrow \mathbf{g}(\mathbf{m}) = \begin{bmatrix} \eta_S \mathbf{g}_S(\mathbf{m}) \\ \eta_P \mathbf{g}_P(\mathbf{m}) \end{bmatrix}. \tag{32}$$

After these modifications, we assume $\boldsymbol{\beta}_i = \beta$ for all $i=1,2,\ldots,M$. Consequently, the Bayesian update formation for $\beta$ can be directly borrowed from [34], in particular,

$$\begin{aligned}\boldsymbol{\beta} &= \dfrac{M}{\mathrm{Tr}\left(\tilde{\mathbf{G}}\boldsymbol{\Sigma}\tilde{\mathbf{G}}^T\right) + \left\|\mathbf{y}^{(n)} - \tilde{\mathbf{G}}\boldsymbol{\mu}\right\|^2} \\ &= \dfrac{M - N + \sum_{i=1}^{N}\boldsymbol{\gamma}_i^{-1}\boldsymbol{\Sigma}_{ii}}{\left\|\mathbf{y}^{(n)} - \tilde{\mathbf{G}}\boldsymbol{\mu}\right\|^2}\end{aligned} \tag{33}$$

Using the determinant identity we can write

$$\log|\mathbf{C}| = -\log\left|\boldsymbol{\Lambda}^{(n+1)}\right| - \log\left|\mathbf{B}^{(n+1)}\right| + \log\left|\boldsymbol{\Lambda}^{(n+1)} + \tilde{\mathbf{G}}^T\mathbf{B}^{(n+1)}\tilde{\mathbf{G}}\right| \tag{34}$$

As a consequence, the gradient of $\log|\mathbf{C}|$ with respect to $\boldsymbol{\gamma}_i$ is

$$\dfrac{\partial \log|\mathbf{C}|}{\partial \boldsymbol{\gamma}_i} = \dfrac{1}{\boldsymbol{\gamma}_i} - \dfrac{\boldsymbol{\Sigma}_{ii}}{\boldsymbol{\gamma}_i^2} \tag{35}$$

Finally, setting the derivative of $\mathcal{L}$ with respect to $\boldsymbol{\gamma}_i$ to zero leads to

$$\dfrac{1}{\boldsymbol{\gamma}_i} - \dfrac{1}{\boldsymbol{\gamma}_i^2}\left(\boldsymbol{\mu}_i^2 + \boldsymbol{\Sigma}_{ii}\right) + \dfrac{\lambda_i}{2} = 0 \tag{36}$$

which readily gives the update formulation

$$\boldsymbol{\Lambda}^{(n+1)} := diag\left\{\gamma_i^{-1}, i=1,2,...,N\right\}, \quad \gamma_i = \begin{cases} \boldsymbol{\mu}_i^2 + \boldsymbol{\Sigma}_{ii} & \text{if } \lambda_i = 0 \\ -\dfrac{1}{2\lambda_i} + \sqrt{\dfrac{1}{4\lambda_i^2} + \dfrac{\boldsymbol{\mu}_i^2 + \boldsymbol{\Sigma}_{ii}}{\lambda_i}} & \text{if } \lambda_i \neq 0 \end{cases} \quad (37)$$

which concludes the derivation of the SBL framework. In the numerical examples of the next section, $\lambda_i = 0$ ($i = 1, 2, ..., N$) are specified. In summary, as listed in Table I and Table II, at ($n$+1)-th iteration of the algorithm, the sensitivity matrix $\mathbf{G}^{(n)}$ and observation $\mathbf{y}^{(n)}$ are estimated based on the previous $\mathbf{m}$, followed by estimating the hyperparameters $\boldsymbol{\Lambda}^{(n+1)}$ and $\mathbf{B}^{(n+1)}$ through (29) or (33) and (37). The coefficients $\boldsymbol{\alpha}^{(n+1)}$ are estimated from (17) and (18).

## IV. NUMERICALL EXPERIMENTS

This paper considers the two-dimensional two-phase (oil/water) immiscible and incompressible multiphase flow, whose differential equations are solved by finite elements method [1-3]. The following simple functions to describe saturation-dependent quantities [1]

$$\lambda_w(s_w) = \frac{(s^*)^2}{\mu_w}, \quad \lambda_o(s_w) = \frac{(1-s^*)^2}{\mu_o}, \quad s^* = \frac{s_w - s_{wc}}{1 - s_{or} - s_{wc}}$$

where $\lambda_w$ and $\lambda_o$ are the water and oil mobility, respectively, $s_{or}$ is the irreducible oil saturation and $s_{wc}$ is the connate water saturation. In this paper, we have assumed $\mu_o = \mu_w = 1.0$ and $s_{or} = s_{wc} = 0$. In addition, just for the purpose of illustrating the proposed algorithms we only consider water injection at the injectors and water/oil production at the producers.

In this section, three sets of synthetic tests resembling waterflooding experiments in oil reservoirs have are carried out to evaluate the performance of the sparse Bayesian algorithm. In each example a 320 m × 320 m × 10 m synthetic reservoir is discretized into a two-dimensional 32 × 32 × 1 uniform grid block system with grid blocks of size 10m × 10 m × 10 m (see Table III). Two-phase (oil/water) immiscible fluid flow simulations are run using a Matlab in-house simulator. A total of 30 observations within uniform intervals of approximately 12 days are considered. The three different well configurations that are used for these experiments are shown in Fig. 2. We refer to the left, middle, and right configurations in Fig. 2 as Reservoir A, B, and C, respectively. Reservoir A (Fig. 3a) portrays a line drive injection using a horizontal well with 32 ports with uniform injection rates (from the left end of the domain). A similar horizontal well with 32 ports is placed at the right end of the domain to produce the displaced oil and water toward the production well. The production ports are under a total production rate constraint to preserve mass balance; these ports produce an equal volume of fluid (oil and water) as the volume of the injected water into the reservoir. Reservoir B (Fig. 3b) includes four injection wells (shown with filled black squares in Fig. 3b) and six production wells (shown with empty black squares). The injectors uniformly inject a total of 1.1 pore volume of water into the reservoir from the left end of the domain during the one year simulation and the producers at the right end of the reservoir uniformly produce a total of 1.1 pore volume of fluids. Reservoir C ( Fig. 3c) includes only one injection wells (shown with filled black squares in Fig. 3c) and five production wells (shown with empty black squares). The injectors uniformly inject a total of one pore volume of water into the reservoir during the one year simulation and the producer uniformly produce a total of one pore volume of fluids. The initial and boundary conditions are assumed to be known perfectly and are listed, along with other important input parameters, in

Table III. For all numerical tests, the initial solution for the permeability is homogeneous with a permeability of 20 mD. We have used the full $32^2$ DCT basis as the compression transform ($\Phi$), in which a sparse solution is sought.

*Numerical Examples 1*: **Reservoir A**

Figure 3 displays the true permeability and the corresponding water saturation profiles after 0, 2, 4, 6, and 12 months. The preferential water flow inside the channels and the resulting early water breakthrough are apparent from these figures. Figure 4a shows the permeability solutions (first row), the corresponding spatial variance (second row), log|DCT| coefficients (third row) and the log(variance) of DCT coefficients at sample iterations of the Bayesian estimation where no sparsity-promoting effect is imposed; that is a Gaussian prior is used for the DCT coefficients. The maximum number of minimization iterations was set to 50, however in most cases after 15 iterations no major improvements in the objective function and the estimated parameters were observed. From Fig.4a, one can find that the reconstruction in this case fails to identify the presence of the channel in the true model. This situation can be further demonstrated from the perspective of solution uncertainty; in particular, the variance results also show that there is very big uncertainty. The distribution of the DCT coefficients in the solution is not sparse as sparsity constraint was not used in this implementation. The estimated saturation profiles in Fig. 5a can not capture the true water front movement and the bypassed oil.

Figure 4b shows sample iterations for the permeability solution (top), the corresponding spatial variance (second row), and the corresponding log|DCT| coefficients (third row) as well as the variance of the estimated DCT coefficients (bottom) for the same inverse problem by using

the proposed Algorithm I, where the sparsity regularization is enforced. In this case, promoting sparsity through sparse priors results in identification of the two channels and their connectivity at the left end of the domain. The results suggest that including sparsity in the solution can lead to better identification of the shape and continuity of the true channels, even though the observations are located at the two ends of the reservoir. A closer look at the variance of the DCT coefficients indicates that although the most significant DCT coefficients are correctly identified a large amount of uncertainty exists in the identified DCT coefficients. The saturation solutions for this experiment are shown in Fig. 5b. The preferential movement of the water front pattern is clearly observed in this case.

Figure 4c shows similar iteration solutions, the corresponding log|DCT| coefficients, and the variance of the DCT coefficients for the same inverse problem that is solved using Algorithm II. A comparison with the results from Algorithm I (Fig. 4b) reveals that a slightly less accurate and more uncertain reconstruction solution is achieved from Algorithm II. The variance in the estimated DCT coefficients is far less in Algorithm II than it is in Algorithm I. Furthermore, while avoiding choosing the artificial parameters as required by Algorithm I, the convergence was achieved more slowly in Algorithm II than it did in Algorithm I. The corresponding saturation profiles for Algorithm II are displayed in Figure 5c. To better appreciate the sparsity-promoting nature of the solutions in Figs. 4b and 4c, Fig. 4d shows the first 400 DCT coefficients of the final solutions. It is evident from Fig. 4d that the significant DCT coefficients have been detected by promoting sparsity through the proposed Algorithms (I and II). Finally, Fig. 6 plots the relative (to initial) pressure reduction for the true and reconstruction results in all

three cases. The match to the observed quantities in the case without sparsity constraint is inferior to those obtained by promoting sparsity.

*Numerical Examples 2*: **Reservoir B**

Here, the reconstruction of the same true model with fewer available measurements is considered. In the new configuration, there are 4 injection wells at the left end of the domain and 6 production wells at the right end (see Fig. 2). Fig. 7 displays the true permeability and the corresponding water saturation profiles after 0, 2, 4, 6, and 12 months. Similar to the results and conclusions in reservoir A can be observed in this example too, the satisfied reconstruction by non-sparse Bayesian estimation cannot be obtained and the corresponding results are not provided here. Fig. 8a shows sample iterations of the permeability, the corresponding log|DCT| coefficients and their variance for algorithm I. While a fast convergence and approximate retrieval of the existing features in the true model is achieved through Algorithm I (i.e. after only 3 iterations), the final solution appears to be blurry with fine scale artifacts, partly because of the non-optimal specification of the parameter $\varepsilon$ in (29). Moreover, the uncertainty in the estimated DCT coefficients is significantly larger in this case, which is consistent with the fact that fewer observations are available for reconstruction. Nonetheless, the overall structure of the true channels seems to be present in the solution.

Figure 8b shows similar results for Algorithm II. In this case the channel structure is correctly captured from limited measurements; however, the solution behaves rather differently in this case. Here, as a result of the decrease in the number of measurements, the solution is mainly focused on very low frequency DCT basis elements to identify only large scale features. The main portion of the estimation variance in this case appears to be concentrated on the low frequency basis elements, implying that estimated DCT coefficients are quite uncertain. A comparison with the true permeability field reveals that the estimated channels are wider than those in the true model and, furthermore, the shape and location of the channels in the middle of the domain where no observations are available are inexact. Nevertheless, the reconstruction results are quite satisfactory given the limited amount of measurements that were used to obtain the solution. The saturation profiles for the three cases discussed in Reservoir B are shown in Fig. 9. As can be confirmed from the figures, the saturation plots are consistent with the permeability estimates that were obtained for this example. Fig. 10 compares the pressure reduction for the true and reconstructed permeability fields for different methods. An noteworthy remark beside the better data match for the case with sparsity constraint is that the match obtained for Algorithm I is better than that of Algorithm II, partly because of the adaptive control of the noise variance of the independent elements of **B**. Before presenting the final example, we note that all the experiments discussed, the only prior information is the sparsity in the DCT domain, which is generally through for any spatially correlated image.

*Numerical Examples 3*: **Reservoir C**

As our final example, we consider a more heterogeneous permeability field under the well configuration in Reservoir C, where only one injection well is located at the upper left upper corner of the investigation domain and five production wells are placed along the opposite edges of the model domain. We only present the results for Algorithm II as the preferred method.

Figure 11a displays the true permeability and the snapshot of the corresponding water saturation profiles. Figure 11b shows sample permeability iterations (top) and the corresponding logarithm of spatial variance (bottom) while figure 11c displays the saturation profiles corresponding to the final solution for this case. It is clear from the reconstruction results that even in the case of more heterogeneous example and when fewer observation locations exist, the major trends in the permeability field are clearly captured by using the sparse Bayesian reconstruction method. While the solution at allocations away from the observations are expected to have large uncertainties in them, the large scale features that can be resolved from the observed quantities are clearly identified by using Algorithm II. The data match for the pressure deduction and saturation measurements at the production wells are plotted in Fig. 11d and 11e, respectively. The overall data match is quite satisfactory except for the pressure reduction at Well P3.

## IV. CONCLUSIONS

We have proposed a sparse Bayesian algorithm for reconstruction of the main features in the rock hydraulic properties from nonlinear dynamic flow measurements. This feature estimation framework is inspired by recognizing that the continuity (correlation) of the spatial distribution in rock properties can be translated into sparse representations in an appropriate compression transform domain (here we used the DCT basis). The transform domain sparse representation summarizes the most salient features in the spatial description of property images and lends itself to effective (feature) estimation techniques that have recently been proposed in sparse reconstruction literature such as the compressed sensing paradigm. While the main theoretical developments in compressed sensing are presented for linear measurements, the underlying principles and the approximate algorithmic developments provide important guidelines for designing sparsity promoting reconstruction techniques in the case of nonlinear (and even dynamic) measurements. Sparsity of the solution in the DCT (or other appropriate transform) domain can be exploited to formulate deterministic regularization methods to improve the solution of ill-posed nonlinear inverse problems. In this paper, we provided a nonlinear version of the probabilistic sparse Bayesian estimation approach for imposing sparsity on the reconstruction solution. This probabilistic framework provide several advantages over the deterministic regularization approach such as elimination of the need to specify the regularization parameter (which can be difficult and computationally expensive to determine) and quantification of the uncertainty in the estimated models and the predictions that are derived from it.

The Laplace distribution which is used as an effective prior for the *$l_1$-norm based sparsity promoting* regularization of the DCT coefficients is indirectly implemented in the sparse Bayesian estimation approach mainly because the Laplace priors are not conjugate to Gaussian likelihoods that are often used to model observations. Therefore, the implemented formulation uses hierarchical priors or Gaussian mixtures to impose sparsity-promoting priors on the DCT coefficients. We extended the application of the sparse Bayesian estimation approach to the nonlinear dynamic problem of characterizing spatial permeability distribution (though we equivalently solved for the DCT domain representation of it) from flow data. We presented alternative formulations for iterative solution of the linearized forms of the original nonlinear problem. Our preliminary results from applying the proposed frameworks to waterflooding examples in two dimensional oil reservoirs suggest that the proposed sparse Bayesian algorithm is indeed and effective solution approach for ill-posed nonlinear inverse problems that are frequently encountered in modeling and identification of subsurface flow and transport systems.

In summary, the results presented in this paper suggest that imposing sparsity constraints via Laplace priors on the solution coefficients (weights) in a Bayesian framework provides a promising approach for estimating spatially correlated subsurface rock hydraulic properties in a compression transform domain. The sparse Bayesian framework is more advantageous to deterministic regularization techniques because it offers a systematic mechanism for include uncertainty in prior and measurements and computing the resulting uncertainty in the estimated parameters and subsequent model predictions. The implications of solving spatial inverse problems in a transform domain using sparsity as a regularization approach or prior information go far beyond characterization of subsurface geologic features (i.e. facies) that we discussed in this paper. We anticipate that the proposed approach can be used to solve several other ill-posed inverse problems that arise in various imaging applications, such as medical imaging, electromagnetic and acoustic inversion, and geophysical tomography where often ill-posed nonlinear inverse problems with nearly sparse unknown parameters are frequently encountered.

**TABLES**

Table I. The procedure for Algorithm I

**Initialization:**
   1. Initialization $\alpha$ and $m$
   2. Selecting suitable sparse transformation $\Phi$
   3. $\lambda_i = 0, i = 1, 2, ..., N$

**While** (stopping criterion not true)
     **Do**

        Computing $y^{(n)}$ and sensitivity matrix $G^{(n)}$

        Finding $\Lambda^{(n+1)} := diag\{\gamma_i^{-1}, i = 1, 2, ..., N\}$

        with $\gamma_i = \mu_i^2 + \Sigma_{ii}$

        Finding $B^{(n+1)} = diag\{\beta_i, i = 1, 2, ..., M\}$

        with $\beta_i = \dfrac{1}{\left(\tilde{G}\Sigma\tilde{G}^T\right)_{ii} + \left(y^{(n)} - \tilde{G}\mu\right)_i^2 + \varepsilon}$

        Computing $\alpha^{n+1} = \Sigma\tilde{G}^T B^{(n+1)} y^{(n)}$ and $\Sigma = \left(\Lambda^{(n+1)} + \tilde{G}^T B^{(n+1)} \tilde{G}\right)^{-1}$

        Computing the convergence criterion
   **END**

Table II. The procedure for Algorithm II

**Initialization:**
    1. Initialization $\alpha$ and $m$
    2. Selecting suitable sparse transformation $\Phi$
    3. $\lambda_i = 0, i = 1, 2, ..., N$

**While** (stopping criterion not true)
    **do**

$$\text{Computing } y = \begin{bmatrix} y_S \\ y_P \end{bmatrix} \leftarrow y = \begin{bmatrix} \eta_S y_S \\ \eta_P y_P \end{bmatrix}$$

$$g(m) = \begin{bmatrix} g_S(m) \\ g_P(m) \end{bmatrix} \leftarrow g(m) = \begin{bmatrix} \eta_S g_S(m) \\ \eta_P g_P(m) \end{bmatrix}$$

and sensitivity matrix $G^{(n)} \leftarrow G^{(n)} = \begin{bmatrix} \eta_S G_S^{(n)} \\ \eta_P G_P^{(n)} \end{bmatrix}$

Finding $\Lambda^{(n+1)} := diag\{\gamma_i^{-1}, i = 1, 2, ..., N\}$

with $\gamma_i = \mu_i^2 + \Sigma_{ii}$

Finding $B^{(n+1)} = diag\{\beta_i, i = 1, 2, ..., M\}$

$$\beta_i = \frac{M - N + \sum_{i=1}^{N} \gamma_i \Sigma_{ii}}{\left\| y^{(n)} - \tilde{G}\mu \right\|^2}, i = 1, 2, ..., M$$

Computing $\alpha^{n+1} = \Sigma \tilde{G}^T B^{(n+1)} y^{(n)}$ and $\Sigma = \left( \Lambda^{(n+1)} + \tilde{G}^T B^{(n+1)} \tilde{G} \right)^{-1}$

Computing the convergence criterion
    **END**

Table III: General simulation information for Reservoir A, B and C

| Parameter | Reservoir A | Reservoir B | Reservoir C |
|---|---|---|---|
| **Simulation Parameters** | | | |
| Phases | Two-phase (o/w) | Two-phase (o/w) | Two-phase (o/w) |
| Simulation Time | 1 years | 1 year | 1 year |
| Grid systems | 32 x 32 x 1 | 32 x 32 x 1 | 32 x 32 x 1 |
| Cell dimensions | 10 x 10 x 10 | 10 x 10 x 10 | 10 x 10 x 10 |
| Rock porosity | 0.20 | 0.20 | 0.20 |
| Initial oil saturation | 1.00 | 1.00 | 1.00 |
| Injection volume | 1.1PV | 1.1PV | 1.0PV |
| Number of injectors | 32 | 4 | 1 |
| Number of producers | 32 | 6 | 5 |
| **Assimilation Information** | | | |
| Observation intervals | 12 days | 12 days | 12 days |
| Obs. at injection wells | Pressure | Pressure | Pressure |
| Obs. at prod. wells | Pressure & saturation | Pressure & saturation | Pressure & saturation |

**FIGURE CAPTIONS**

**Figure 1.**

**(a)** The reservoir permeability distribution, its log-|DWT| coefficients, and the reconstructions with 1%, 2% and 5% largest DCT coefficients.

**(b)** The reservoir permeability distribution, its log-|DCT| coefficients, and the reconstructions with 1%, 2% and 5% largest DCT coefficients.

**Figure 2.** The sketch map of two considered well configuration

(a) The sketch map of Reservoir A

(b) The sketch map of Reservoir B, where four injection wells shown with filled black squares and six production wells shown with empty black squares.

(c) The sketch map of Reservoir C, where one injection wells shown with filled black squares and five production wells shown with empty black squares.

**Figure 3.**

The true permeability and the saturation after 0months, 2 months, 4 months, 6 months and 12 months for Reservoir A.

**Figure 4.**

(a) The reconstructed permeability of reservoir A using Bayesian estimation without sparsity constraint after 1 iteration, 3 iterations, 6 iterations, 9 iterations, 12 iterations and 15 iterations (top), the corresponding spatial variance (second row), log|DCT| coefficients (third row) and the log of variance of DCT coefficients(bottom).

(b) The reconstructed permeability of reservoir A with Algorithm I after 1 iteration, 3 iterations, 6 iterations, and 9 iterations, the corresponding the spatial variance (second row), the distributions of log|DCT| coefficients (third row) and the variance of DCT-coefficents (fourth row).

(c) The reconstructed permeability of reservoir A with Algorithm II after 1 iteration, 9 iterations, 12 iterations, 21 iterations, 27 iterations and 33 iterations, the corresponding the spatial variance

(second row), the distributions of log|DCT| coefficients (third row) and the variance of DCT-coefficents (fourth row).

(d) The first row shows the distribution of log|DCT| of true, and reconstructions with traditional Bayesian estimation without sparsity constraint, Algorithm I and Algorithm II; the left figure of second row shows the comparisons of true DCT coefficients (with black line) and reconstructed coefficients (with red line) using Algorithm I; the right figure of second row shows the comparisons of true DCT coefficients (with black line) and reconstructed coefficients (with red line) using Algorithm II.

**Figure 5.**

(a)The reconstructed permeability and its corresponding saturation after 0 months, 2 months, 4 months 6 months and 12 months for traditional Bayesian estimation without sparsity constraint

(b)The reconstructed permeability and its corresponding saturation after 0 months, 2 months, 4 months 6 months and 12 months for Algorithm I

(c)The reconstructed permeability and its corresponding saturation after 0 months, 2 months, 4 months 6 months and 12 months for Algorithm II.

**Figure 6.**

(a)The comparison of pressure reduction at production well #1 between true one and reconstructions with different algorithm.

(b)The comparison of pressure reduction at production well #8 between true one and reconstructions with different algorithm.

(c)The comparison of pressure reduction at production well #24 between true one and reconstructions with different algorithm.

(d)The comparison of pressure reduction at production well #32 between true one and reconstructions with different algorithm.

**Figure 7.**

The true permeability and the saturation after 0months, 2 months, 4 months, 6 months and 12 months for Reservoir B.

**Figure 8.**

(a)The reconstructed permeability of reservoir A with Algorithm I after 1 iteration, 3 iterations, 6 iterations, 9 iterations, 12 iterations, and 15 iterations, the corresponding the spatial variance (second row), the distributions of log|DCT| coefficients (third row) and the variance of DCT-coefficients (fourth row).

(b)The reconstructed permeability of reservoir A with Algorithm II after 1 iteration, 9 iterations, 12 iterations, 34 iterations, 45 iterations and 48 iterations, the corresponding the spatial variance (second row), the distributions of log|DCT| coefficients (third row) and the variance of DCT-coefficients (fourth row).

(c)The first row shows the distribution of log|DCT| of true, and reconstructions with Algorithm I and Algorithm II; the left figure of second row shows the comparisons of true DCT coefficients (with black line) and reconstructed coefficients (with red line) using Algorithm I; the right figure of second row shows the comparisons of true DCT coefficients (with black line) and reconstructed coefficients (with red line) using Algorithm II.

**Figure 9.**

(a)The reconstructed permeability and its corresponding saturation after 0 months, 2 months, 4 months 6 months and 12 months for Algorithm I

(b)The reconstructed permeability and its corresponding saturation after 0 months, 2 months, 4 months 6 months and 12 months for Algorithm II.

**Figure 10.**

(a)The comparison of pressure reduction at production well #1 between true one and reconstructions with different algorithms.

(b)The comparison of pressure reduction at production well #2 between true one and reconstructions with different algorithms.

(c)The comparison of pressure reduction at production well #3 between true one and reconstructions with different algorithms.

(d)The comparison of pressure reduction at production well #4 between true one and reconstructions with different algorithms.

**Figure 11.**

(a) The true permeability and the saturation after 0months, 2 months, 4 months, 6 months and 12 months for Reservoir B.

(b) It shows the distribution of the iteration samples (first row) and corresponding log of variance (second row) by using Algorithm II.

(c) The reconstructed permeability and its corresponding saturation after 0 months, 2 months, 4 months 6 months and 12 months.

(d)The comparison of pressure reduction at production wells between true one (solid line) and reconstruction (dashed line).

(e)The comparison of saturation at production wells between true one (solid line) and reconstruction (dashed line).

**FIGURES**

**Fig.1**

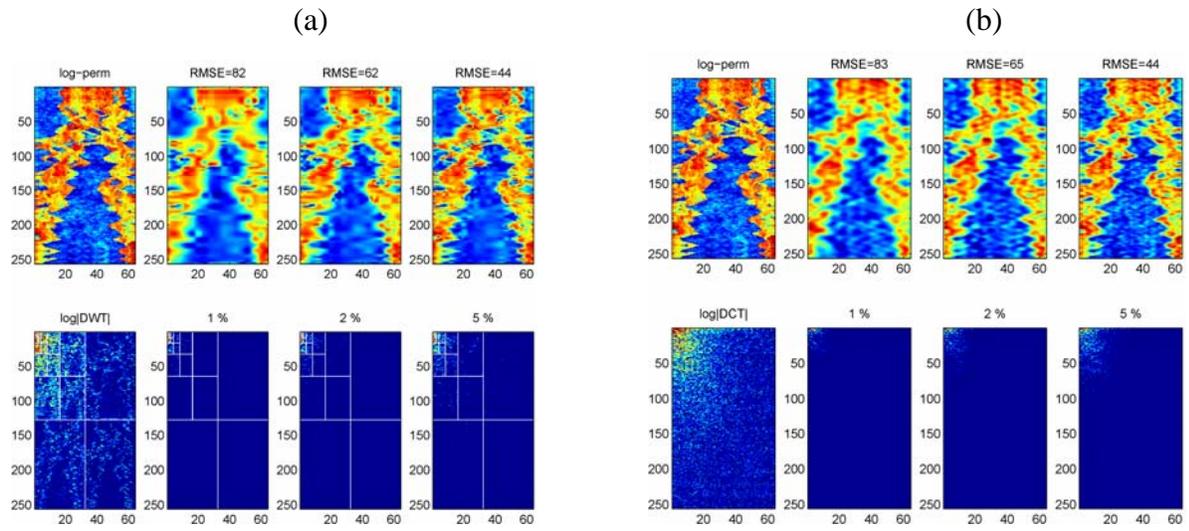

**Fig.2**

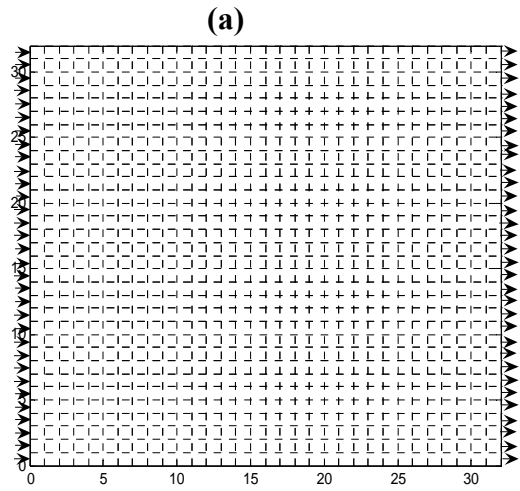

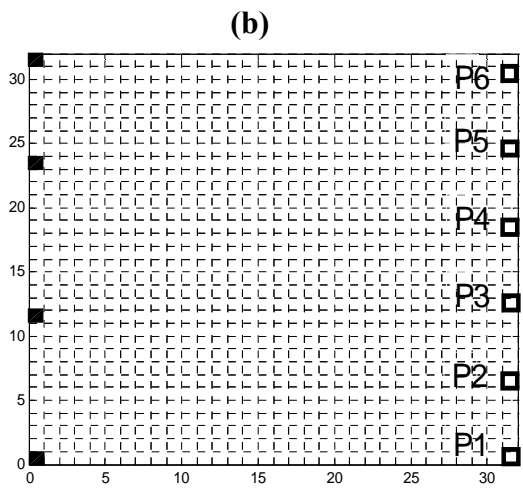

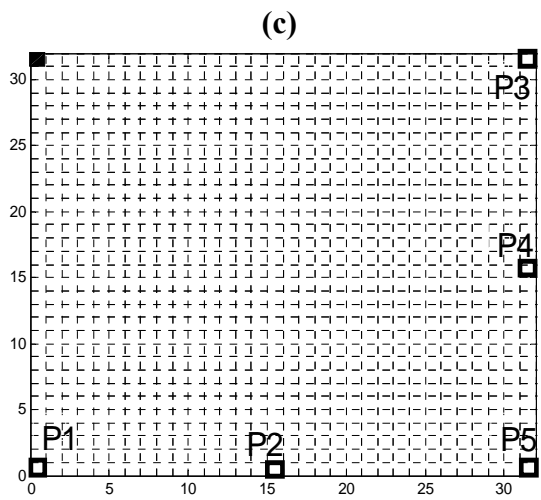

**Fig.3**

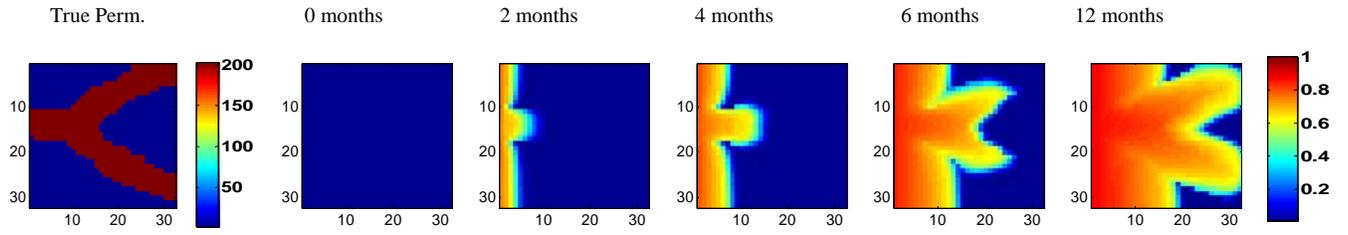

**Fig. 4**

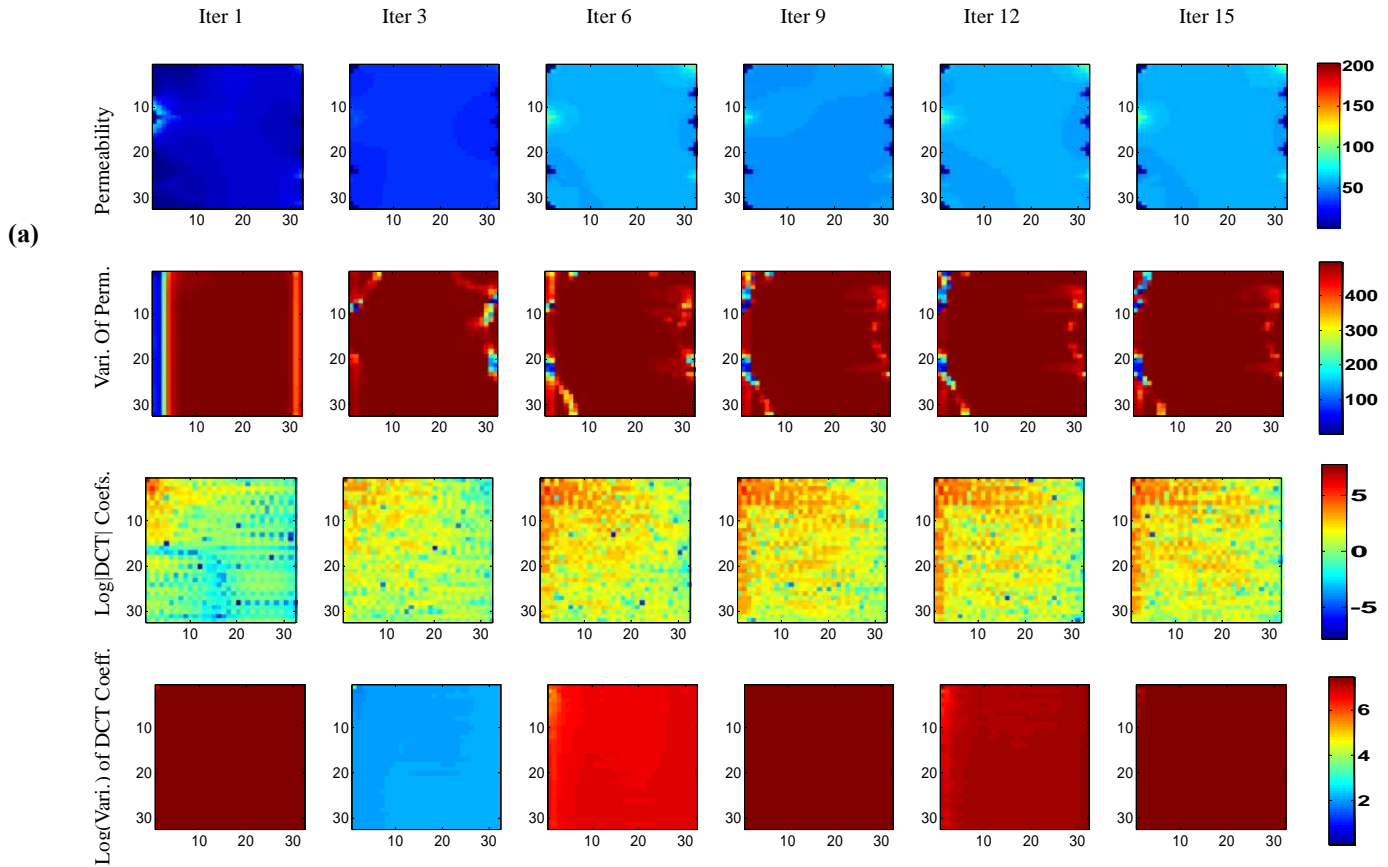

(a)

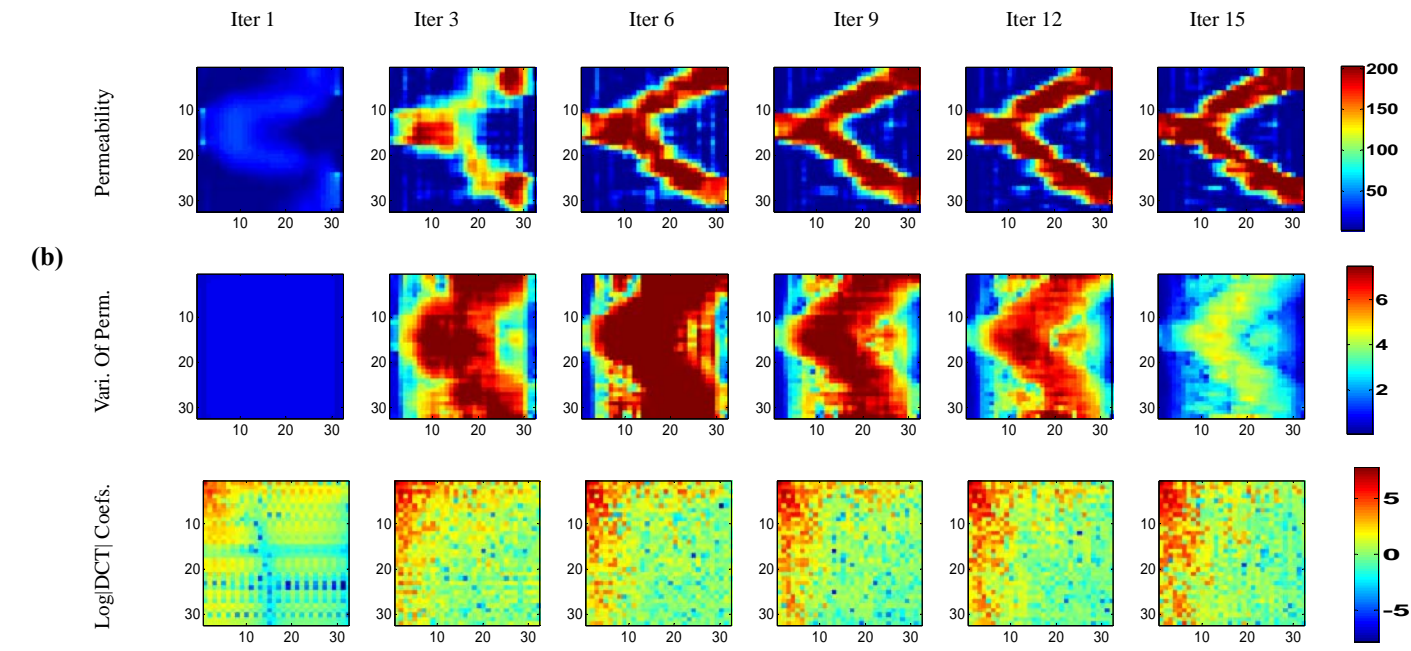

(b)

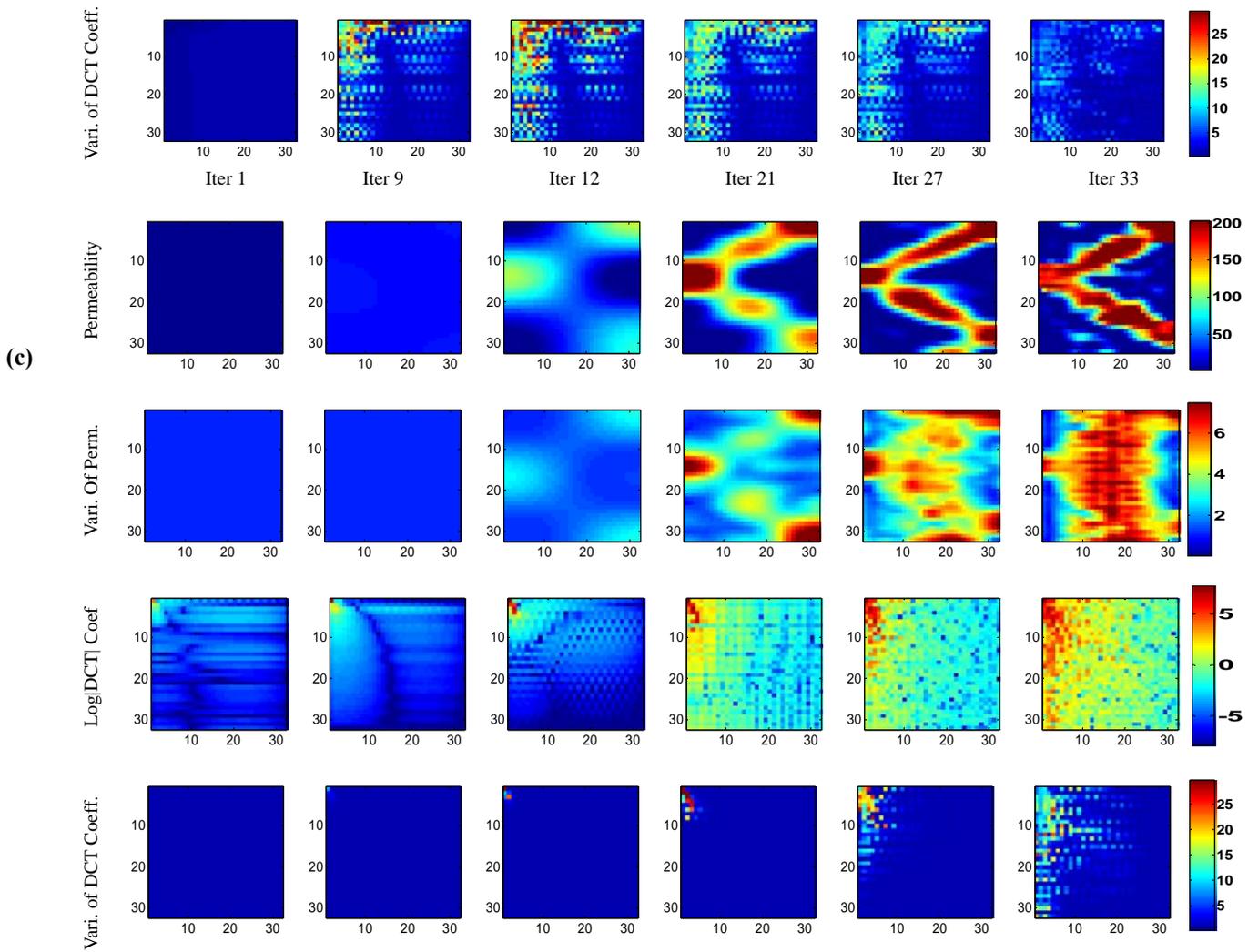

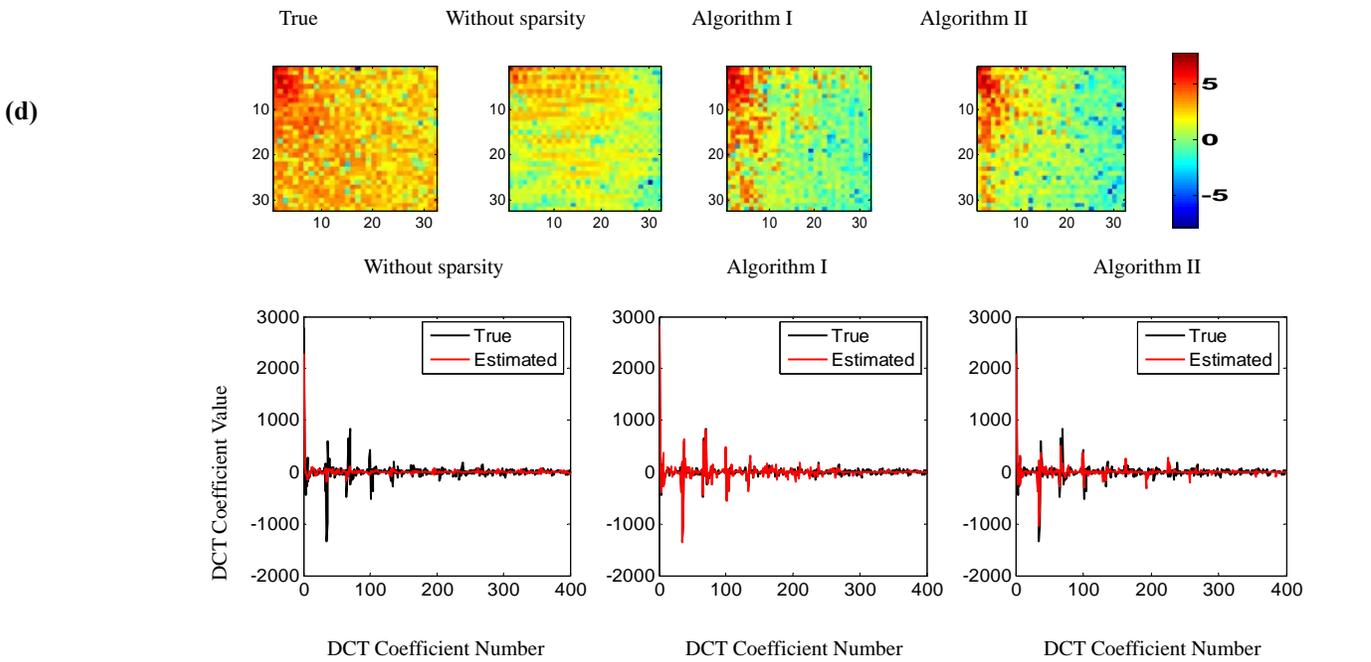

**Fig.5**

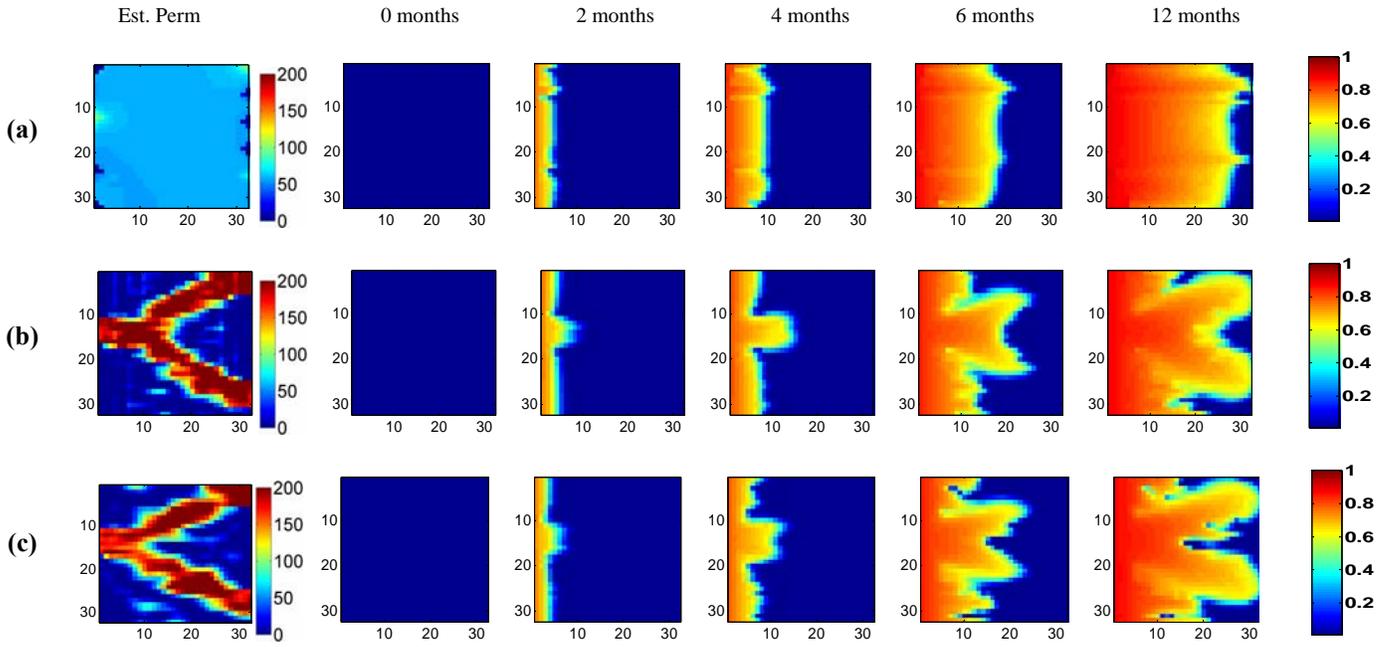

**Fig. 6**

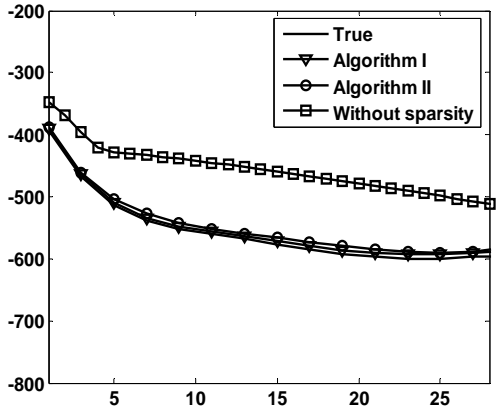

(a) Well P1

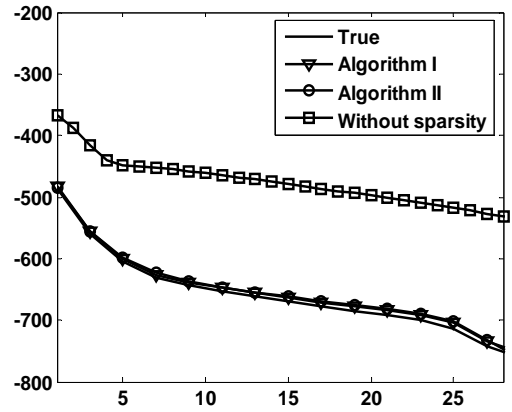

(b) Well P8

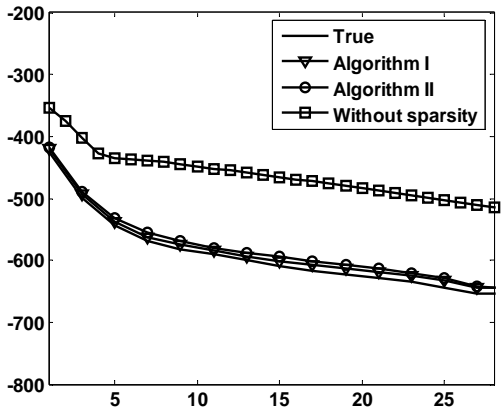

(c) Well P24

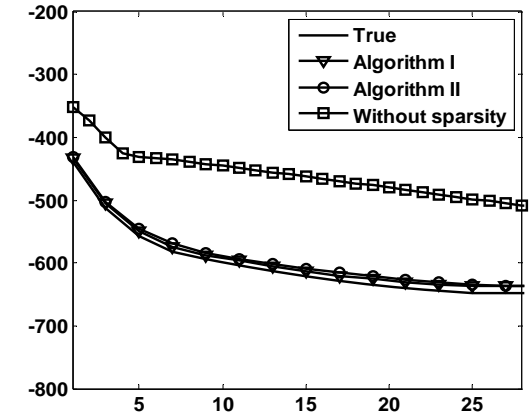

(d) Well P32

**Fig.7**

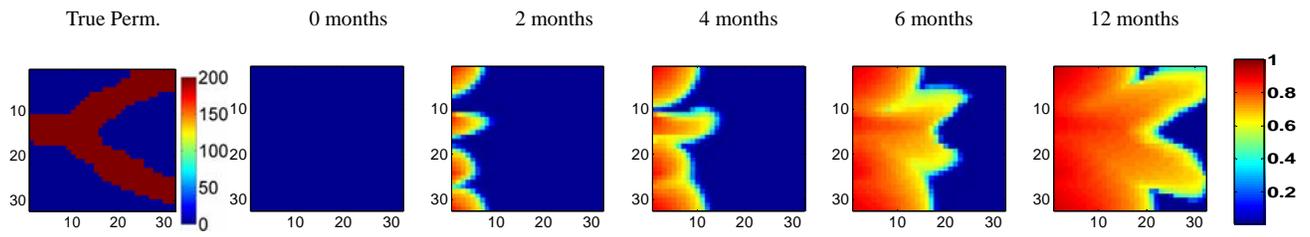

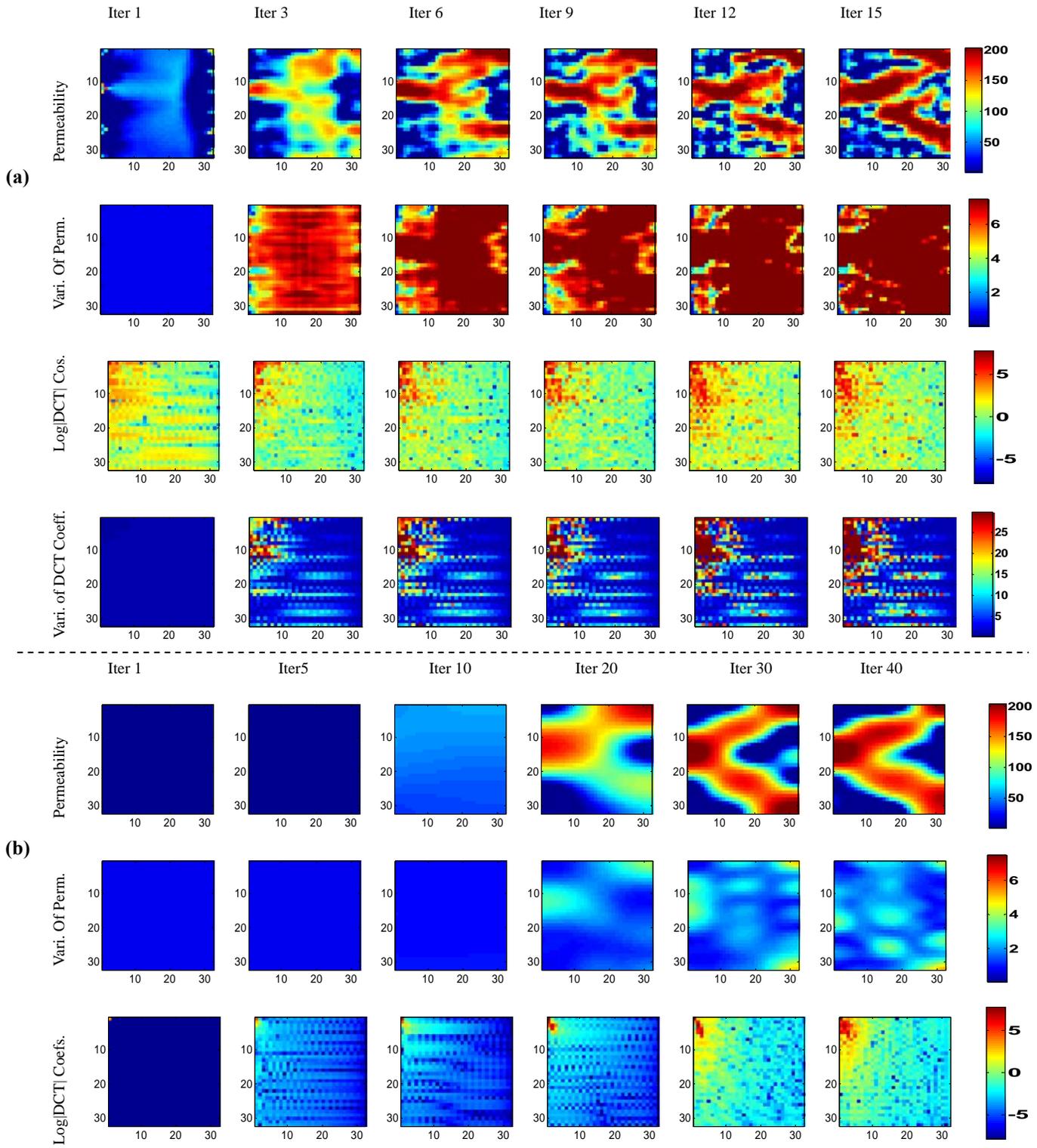

Fig. 8

**(c)**

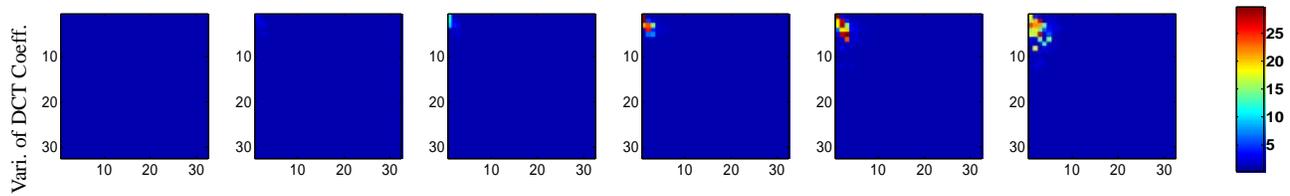
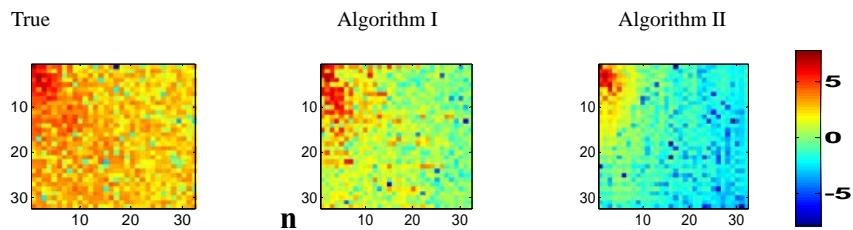
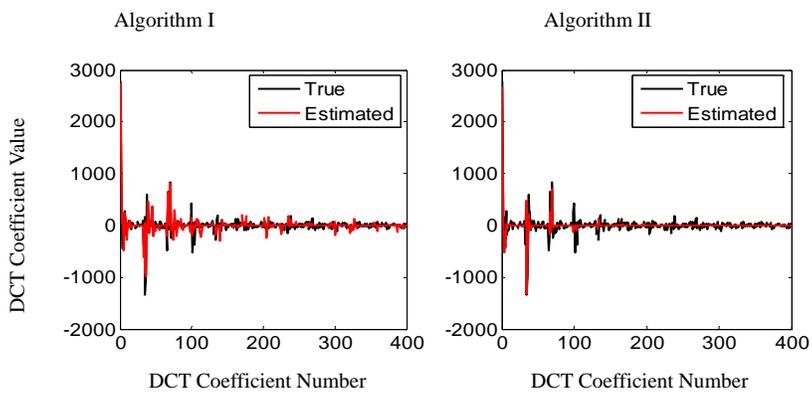

**Fig.9**

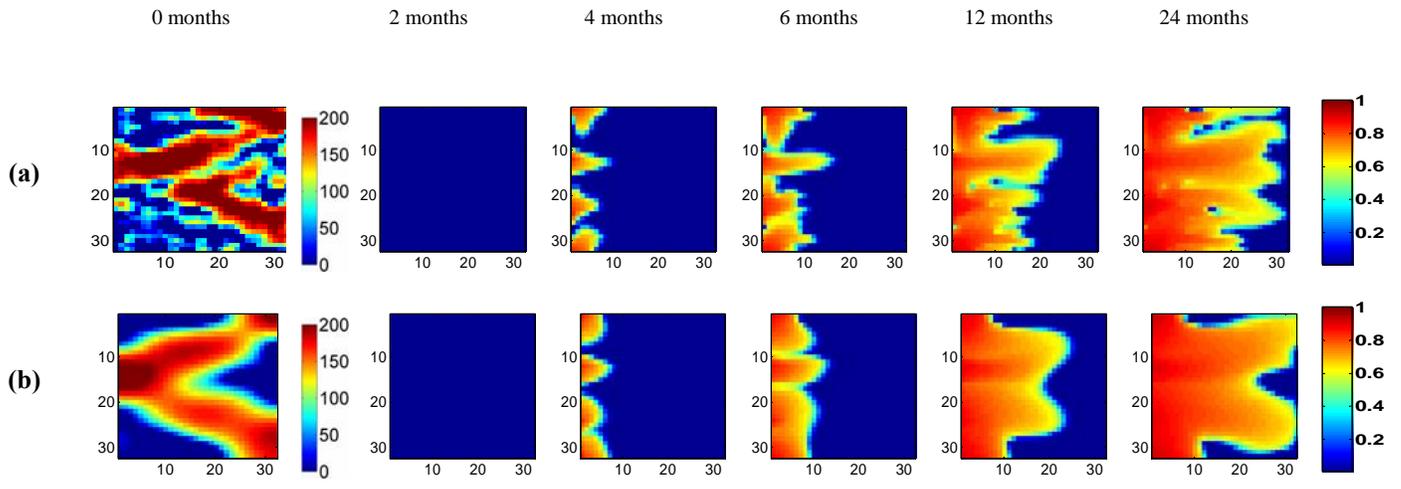

(a)

(b)

0 months    2 months    4 months    6 months    12 months    24 months

**Fig.10**

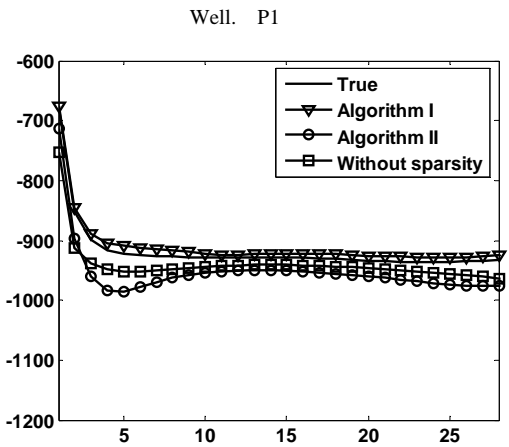
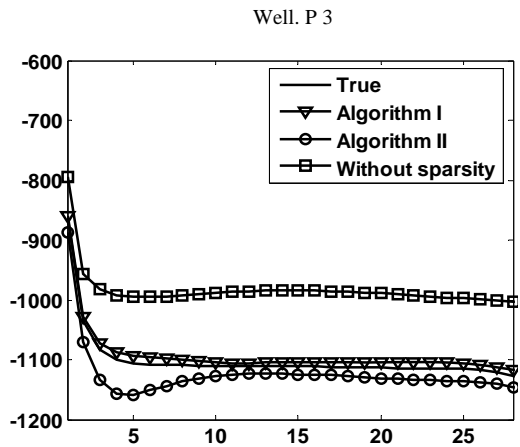
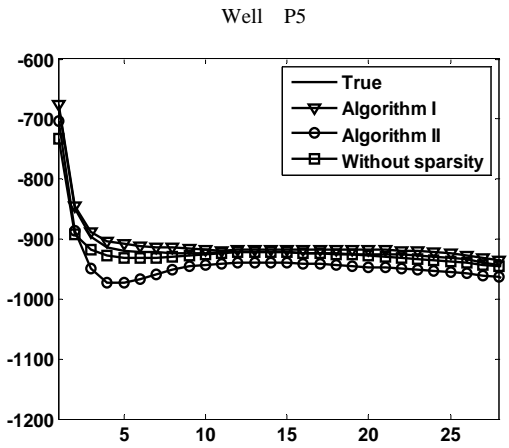
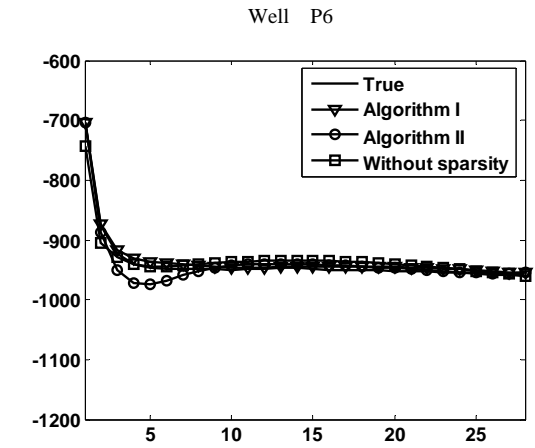

**Fig.11**

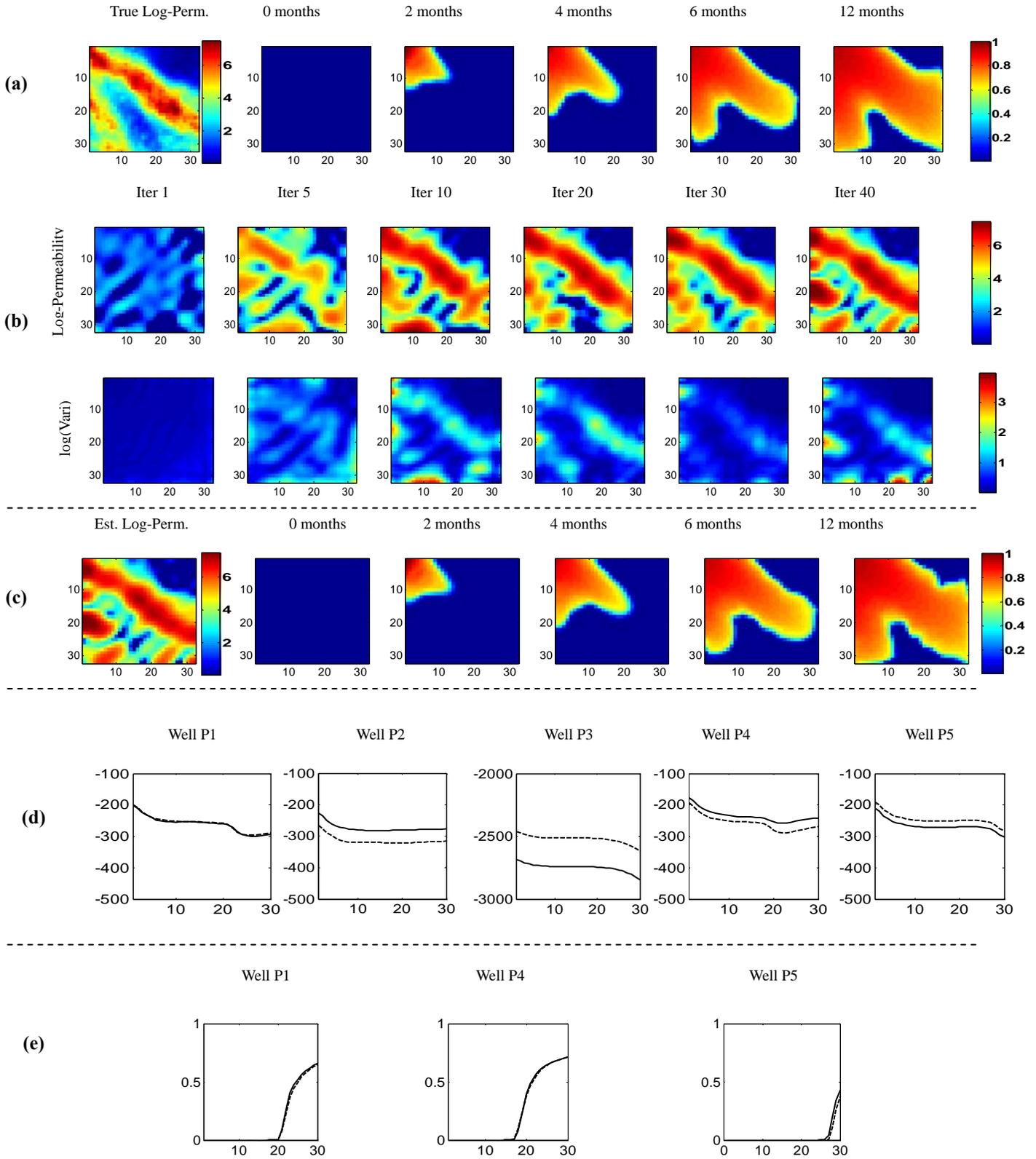